\documentclass{amsart}\usepackage[reset,a4paper,vmargin=5truecm,hmargin=4truecm]{geometry}
\usepackage{amsmath,amsthm,amssymb}

\usepackage{graphicx}
\usepackage{url}
\usepackage{cite}

\newtheorem{thm}{Theorem}
\newtheorem{lem}{Lemma}
\newtheorem{cor}{Corollary}

\newtheorem*{rem}{Remark}

\begin{document}

\title[Tiling Problem]{Tiling Problem: Convex Pentagons for Edge-to-Edge 
Monohedral Tiling and Convex Polygons for Aperiodic Tiling}

\address{The Interdisciplinary Institute of Science, Technology and Art\\
Suzukidaini-building 211, 2-5-28 Kitahara, Asaka-shi, Saitama, 351-0036, 
Japan}
\email{ismsugi@gmail.com}
\author{Teruhisa Sugimoto}
\maketitle	

\begin{abstract}
We show that convex pentagons that can generate edge-to-edge 
monohedral tilings of the plane can be classified into exactly eight 
types. Using these results, it is also proved that no single convex 
polygon can be an aperiodic prototile without matching conditions 
other than ``edge-to-edge.''
\end{abstract}

\section{Introduction}

A collection of sets (the ``tiles'') is a \textit{tiling} of the plane if their union 
is the entire plane, but the interiors of different tiles are disjoint. If all 
tiles in a tiling are of the same size and shape, then the tiling is called 
\textit{monohedral}, and the polygon in the monohedral tiling is called 
the \textit{prototile} of the monohedral tiling, or simply, the 
\textit{polygonal tile}~\cite{G_and_S_1987, Sugimoto_2012a, Sugimoto_2012b, 
Sugimoto_2015, Sugimoto_III}. 

In the problem of classifying the types of convex pentagonal tiles, only the 
pentagonal case remains unsettled. At present, known convex pentagonal tiles 
are classified into essentially 15 different types (see type 1, type 2, ... , type 15 
in~\cite{G_and_S_1987, Mann_2015, Sugimoto_2012a, Sugimoto_2015, 
Sugi_Ogawa_2006, wiki_pentagon_tiling}). 
However, it is not known whether this list of types is complete (perfect) or 
not\footnote{In 1918, Reinhardt showed five types of convex pentagonal 
tiles in his dissertation~\cite{Reinhardt_1918}. In 1968, Kershner published 
the article which showed eight types of convex pentagonal tiles which 
added three new types~\cite{Kershner_1968}.  For several years, it was 
generally accepted as fact that Kershner's list of types of convex pentagonal 
tiles is complete~\cite{G_and_S_1987}. However, after Gardner wrote an 
expository article on convex polygonal tiles in 1975~\cite{Gardner_1975a}, 
Richard James discovered one type and Marjorie Rice discovered 
four types~\cite{Gardner_1975b, G_and_S_1987, Schatt_1978}. In 1985, the 
14th type was discovered by Rolf Stein~\cite{G_and_S_1987, Stein_1985,Wells_1991}. 
In July, 2015, Casey Mann et al. discovered the 15th type~\cite{Mann_2015, 
wiki_pentagon_tiling}. See Appendix \ref{appA1} for details of type 15.}~\cite{Bagina_2004, 
G_and_S_1987, Hallard_1991, Hirshh_1985, Schatt_1978, 
Sugimoto_2012a, Sugimoto_2012b, Sugimoto_2015, Sugimoto_III, Sugi_Ogawa_2006, 
Wells_1991, wiki_pentagon_tiling}. Note that the problem 
of classifying the types of convex pentagonal tiles and the problem of classifying the 
pentagonal tilings (tiling patterns) are quite different. 

Tiling by convex polygons is called \textit{edge-to-edge} if any two polygons 
are either disjoint or share one vertex or one entire edge in common. In an 
edge-to-edge tiling, the vertices are called \textit{nodes}. Thus, the 
vertices of several polygons meet at a node of an edge-to-edge tiling, 
and the number of polygons meeting at the node is called the 
\textit{valence} of the node~\cite{Sugimoto_2012a, Sugimoto_2012b, Sugimoto_2015, 
Sugimoto_III, Sugi_Ogawa_2006}.

In tilings by the 15 types of convex pentagonal tiles, there are 
edge-to-edge and non-edge-to-edge tilings. A convex pentagonal tile 
belonging only to type 3, type 10, type 11, type 12, type 13, type 14, or type 15 
cannot generate an edge-to-edge tiling. The other eight types can generate 
an edge-to-edge tiling. (In Table 1, the representative tilings of type 1 or type 2 
are non-edge-to-edge as shown, but there exist such convex pentagons that can 
generate edge-to-edge tilings belonging to type 1 or type 2~\cite{Sugimoto_2012a, 
Sugimoto_2015, Sugi_Ogawa_2006}). Our main result is the following theorem\footnote{ 
We know that the same result was obtained by Bagina in 2011 after this result 
was obtained in 2012~\cite{Bagina_2011, Sugimoto_2012b}.}~\cite{Sugimoto_2012b, 
Sugimoto_III}.

\begin{thm}\label{t1}
If a convex pentagon can generate an edge-to-edge monohedral tiling, 
then it belongs to one of eight types in Table$~\ref{EEpenta}$.
\end{thm}

\begin{rem}
These eight types are not necessarily ``disjoint.'' Some of 
them can contain the same convex pentagon.
\end{rem}

\begin{table}[p]
{\small
\caption{Eight types of convex pentagonal tiles that can generate 
edge-to-edge tilings.} 
\label{EEpenta}
}
\centering
  \begin{tabular}{c}  
    \begin{minipage}{12.8cm}
      \centering\scalebox{0.73}{\includegraphics{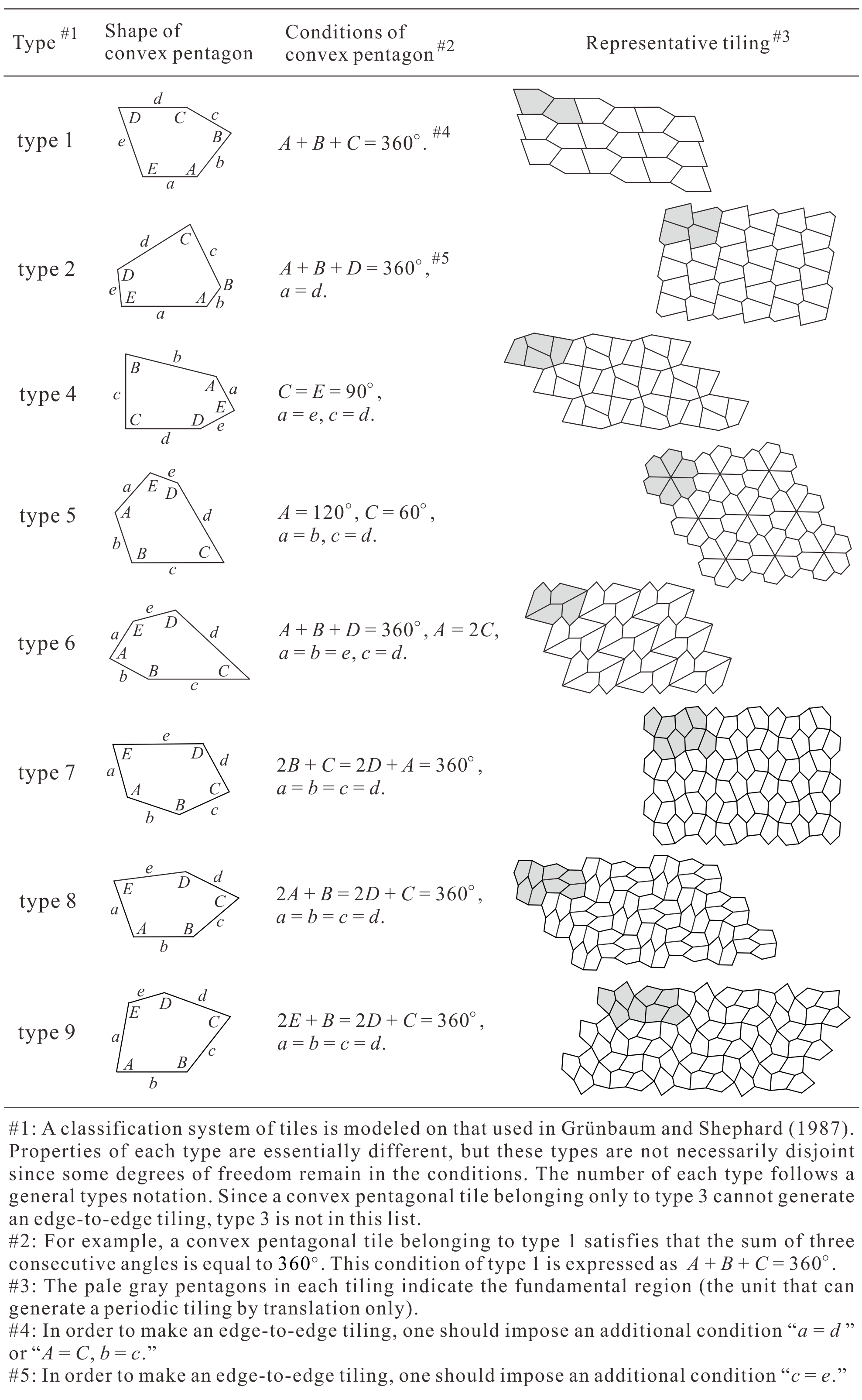}}
    \end{minipage} 
  \end{tabular}
\end{table}

\medskip
A tiling of the plane is \textit{periodic} if the tiling can be translated 
onto itself in two nonparallel directions. More precisely, a tiling is periodic 
if it coincides with its translation by a nonzero vector. A set of prototiles 
is called \textit{aperiodic} if congruent copies of the prototiles admit infinitely 
many tilings of the plane, none of which are periodic. Thus, an aperiodic set 
of prototiles admits only nonperiodic tilings. (A tiling that has no periodicity is 
called nonperiodic. On the other hand, a tiling by aperiodic sets of prototiles 
is called aperiodic. Note that, although an aperiodic tiling is a nonperiodic 
tiling, a nonperiodic tiling is not necessarily an aperiodic tiling.) For 
example, the Penrose tiling is a nonperiodic tiling, and it can also be 
considered an aperiodic tiling that is generated by the aperiodic set of 
prototiles with a matching condition~\cite{Akiyama_2012, G_and_S_1987, Hallard_1991, 
Socolar_2011}. (A matching condition is a condition concerning 
edge-matching. In some cases, it can be represented by assigning colors 
and orientations to several edges of the prototiles.) For aperiodic tilings, 
there are the following two problems~\cite{Hallard_1991}.

\medskip
\noindent
(i) Is there a single aperiodic prototile (with or without a matching 
condition), that is, one that admits only aperiodic tilings by congruent 
copies?

\noindent
(ii) It is well known that there is a set of three convex polygons that are 
aperiodic with no matching condition on the edges. Is there a set of 
prototiles with a size less than or equal to two that is aperiodic?

\medskip
\noindent
For Problem (i), the affirmative solution (a convex hexagon with a matching 
condition) was proved in recent years~\cite{Socolar_2011}. The aperiodic tiling by 
a set of three convex polygons (one hexagon and two pentagons) with no matching 
condition on the edges in Problem (ii) is based on the Penrose tiling~\cite{G_and_S_1987}. 
Note that the aperiodic tiling is edge-to-edge.

Using Theorem \ref{t1}, it can be proved that every convex polygon that generates 
an edge-to-edge tiling can generate a periodic tiling. Hence, we have the 
following theorem\footnote{ Any triangle and any convex quadrilateral can 
generate a periodic edge-to-edge tiling. All convex hexagons that can 
generate a monohedral tiling are categorized into three types and can 
generate a periodic edge-to-edge tiling. There is no convex polygon with 
seven or more edges that can generate a monohedral tiling~\cite{Gardner_1975a, 
G_and_S_1987, Hallard_1991, Klamkin_1980, Reinhardt_1918}.}~\cite{Sugimoto_2012b}.

\begin{thm}\label{t2}
Without matching conditions other than ``edge-to-edge,'' 
no single convex polygon can be an aperiodic prototile.
\end{thm}

\section{Outline of the Proof of Theorem 1}

Our proof is based on Bagina's Proposition (see Proposition 2.1 in 
\cite{Bagina_2004}, Section 1 in \cite{Sugimoto_2012a} or,  Subsection 
2.2 in \cite{Sugimoto_2015}), which implies the following: 

\noindent
\textit{In any edge-to-edge tiling of the plane by convex pentagonal 
tiles, there is a pentagon around which at least three nodes have 
valence $3$; in other words, there exists a tile with at least three 
$3$-valent nodes}.

Now, let $G$ be a convex pentagonal tile candidate that can generate an 
edge-to-edge tiling. Then, by Bagina's Proposition, it has at least three 
vertices that will become 3-valent nodes in the tiling. We choose (any) two 
of them as $V_{1}$ and $V_{2}$ and let $v_{1}$ and $v_{2}$ be the conditions on 
the angles at $V_{1}$ and $V_{2}$, respectively. (For example, if $V_1 = A$, 
then $v_{1}$ could be $A + B + D = 360^ \circ $, $2A + B = 360^ \circ $, 
etc.) By substituting all possible angle conditions in $v_{1}$ and $v_{2}$ and 
considering the lengths of the matching edges, we can produce 465 patterns 
of pentagons as the candidates $G$. Examining these pentagons one by one, 
we sort them into (i) geometrically impossible cases, (ii) the cases that 
cannot generate an edge-to-edge tiling, (iii) the eight types in Table 1, 
and (iv) uncertain cases. In this first-stage sorting, there remain 34 
uncertain cases (for details, see \cite{Sugimoto_2012a}). These remaining 34 
cases are further refined into 42 cases by imposing extra edge 
conditions~\cite{Sugimoto_2015}. To check these 42 cases, we further 
prepare the following lemma and corollary~\cite{Sugimoto_2015, Sugimoto_III}.

\begin{lem}\label{lem1}
If the density of k-valent nodes for $k \ge 5$ in an edge-to-edge 
monohedral tiling by a convex pentagon is greater than zero, then 
there exists a tile with four or more $3$-valent nodes.
\end{lem}

\begin{cor}\label{cor1}
If a convex pentagon with exactly three vertices that can 
simultaneously belong to tentative $3$-valent nodes is a 
convex pentagonal tile that can generate an edge-to-edge tiling, 
then an edge-to-edge monohedral tiling by the tile is a tiling of 
the plane with only $3$- and $4$-valent nodes.
\end{cor}

\medskip
\noindent
Note that the density of $k$-valent nodes is the ratio of the number of 
$k$-valent nodes to the number of pentagons in an edge-to-edge monohedral 
tiling by a convex pentagon (see Subsection 2.1 of  \cite{Sugimoto_2015} for the 
exact definition of the density).

Applying Bagina's Proposition, Lemma \ref{lem1}, Corollary \ref{cor1}, etc., we 
could conclude that every convex pentagon in the remaining 42 cases either 
belongs to the eight types in Table 1 or cannot generate an edge-to-edge tiling. 
A computer was used to perform such an investigation.

\bigskip
\noindent
{\textbf{Acknowledgments.}
The author would like to thank Emeritus Professor Hiroshi 
Maehara, University of the Ryukyus, for his helpful comments and 
suggestions, Professor Shigeki Akiyama, 
University of Tsukuba, for his mathematical helpful comments 
and guidance, and Emeritus Professor Masaharu Tanemura, 
The Institute of Statistical Mathematics, for a critical 
reading of the manuscript.

\bigskip
\appendix
\def\thesection{\Alph{section}}
\section{}
\label{appA1}

In July, 2015, Casey Mann, Jennifer McLoud-Mann, and David Von Derau 
discovered the convex pentagonal tile of type 15 (see Figure~\ref{type15fig}) 
by a computerized search~\cite{Mann_2015, wiki_pentagon_tiling}.

\renewcommand{\figurename}{{\small Figure.~}}
\begin{figure}[!h]
\centering \includegraphics[width=12cm,clip]{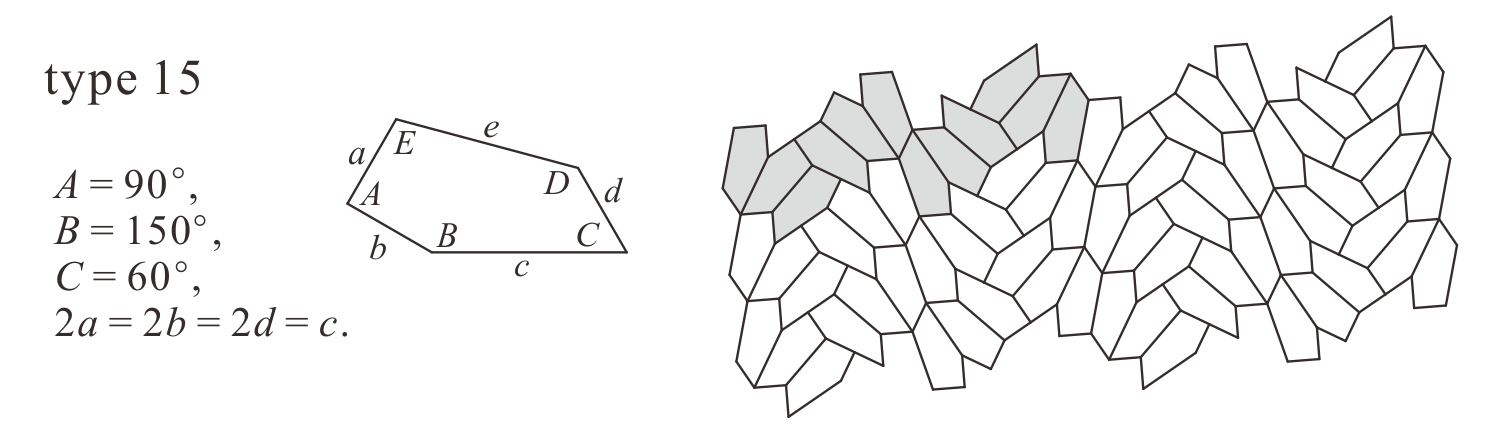}
\caption{{\small 
Convex pentagonal tiles of type 15.} 
\label{type15fig}
}
\end{figure}

\end{document}